\documentclass[draft]{agujournal2018}

\usepackage[inline]{trackchanges}

\usepackage{apacite}
\usepackage{url} 
\usepackage{lineno}
\usepackage{cancel}

\usepackage{amsfonts} 
\usepackage{amsmath} 

\newcommand{\bo}[1]{{\bf #1}}

\journalname{Submitted to Journal of Advances in Modeling Earth Systems (JAMES)}

\begin{document}

\title{A method for generating coherent spatially explicit maps of seasonal palaeoclimates from site-based reconstructions}

%
%



\authors{Cleator, S.F.\affil{1}, Harrison, S.P.\affil{2}, Nichols, N.K.\affil{3}, \\ Prentice, I.C.\affil{4} and Roulstone, I.\affil{1}}

\affiliation{1}{Department of Mathematics, University of Surrey, Guildford GU2 7XH, UK}
\affiliation{2}{School of Archaeology, Geography and Environmental Sciences (SAGES), University of Reading, Whiteknights, Reading RG6 6AH, UK}
\affiliation{3}{Department of Mathematics \& Statistics, University of Reading, Whiteknights, Reading RG6 6AX, UK}
\affiliation{4}{AXA Chair of Biosphere and Climate Impacts, Department of Life Sciences, Imperial College London, Silwood Park Campus, Buckhurst Road, Ascot SL5 7PY, UK}





\correspondingauthor{Sean F. Cleator}{s.cleator@surrey.ac.uk}



\begin{keypoints}
\item 3D-Variational technique used to reconstruct seasonal palaeoclimates
\item Results reinforce spatial and seasonal coherence of temperature and moisture anomalies
\item Applied to southern Europe at the Last Glacial Maximum, ca 21000 years ago
\end{keypoints}

%
%

\begin{abstract}
    We describe a method for reconstructing spatially explicit maps of seasonal palaeoclimate variables from site-based reconstructions.
    Using a 3D-Variational technique, the method finds the best statistically unbiased, and spatially continuous, estimate of the palaeoclimate anomalies through combining the site-based reconstructions and a prior estimate of the palaeoclimate state.
    By assuming a set of  correlations in the error of the prior, the resulting climate is smoothed both from month to month and from grid cell to grid cell.
    The amount of smoothing can be controlled through the choice of two length-scale values.
    The method is applied to a set of reconstructions of the climate of the Last Glacial Maximum (ca. 21,000 years ago, yr BP) for southern Europe derived from pollen data with a prior derived from results from the third phase of the Palaeoclimate Intercomparison Project (PMIP3). 
    We demonstrate how to choose suitable values for the smoothing length scales from the datasets used in the reconstruction. 
\end{abstract}

%
%

\section{Introduction}

Past climates provide useful examples of how the climate system has responded to changes in external forcing, such as orbitally-induced changes in incoming solar radiation, and internal Earth system feedbacks, such as changes in atmospheric CO$_2$ concentration ([CO$_2$]) or ice sheet extent \citep{harrison2012records}.
Reconstructions of past climate states are now routinely used to evaluate the performance of the climate models that are used to project the trajectory of future climate changes \citep{harrison2014climate, schmidt2014using, harrison2015evaluation, kageyama2018pmip4}. 
The Last Glacial Maximum (LGM, ca. 21,000 years ago) has been a major focus for these evaluations because the change in climate forcing was as large (albeit different in type) as "high-end" changes projected for the end of the 21st century \citep{braconnot2012evaluation, kageyama2018pmip4}. 
These evaluations obviously depend on the availability of quantitative reconstructions of key climate variables and this has led to the creation of benchmark data sets documenting climate conditions over land (e.g. \citealp{bartlein2011pollen}) and ocean (e.g. \citealp{margo2009constraints}).

Past climate conditions can be inferred from environmental records which respond to climate, including sedimentological, geomorphological, chemical, isotopic and biological records \citep{gornitz2008encyclopedia, bradley1999paleoclimatology}. 
Quantitative reconstructions of climate variables can be obtained from these records either using statistical techniques based on modern day climate-response relationships (e.g. \citealp{braak1993weighted}; see also discussion in \citealp{bartlein2011pollen}) or by inversion of a model that simulates the response of a particular type of environmental record to climate (e.g. \citealp{garreta2009method, steiger2017climate}). 
Pollen preserved in anoxic lake and bog sediments through time is the most widespread source of data for the reconstruction of terrestrial climates \citep{bartlein2011pollen, marsicek2018reconciling}, because pollen abundance reflects the distribution of different plant taxa that have highly specific climatic requirements \citep{woodward1987climate, harrison2010ecophysiological} and the pollen-preserving sediments can be accurately dated using radiometric techniques.
One important characteristic of all of the environmental records that are used for climate reconstruction, including pollen, is that both the primary data and the climate reconstructions are generated for individual sites. 
Geological and climatic factors mean that the distribution of potential sites is spatially non-uniform: speleothem records, for example, are confined to karst areas; pollen preservation requires anoxic conditions and thus pollen records are not common in arid regions.
Furthermore, issues of accessibility and scientific interests means that the actual sampling of potential sites is non-uniform, so there are often large geographic gaps in the data coverage \citep{harrison2016what}.
While pollen records, for example, are well-sampled across Europe and North America, there are far fewer records from central Eurasia or the tropics.
Furthermore, geological preservation issues mean that the number of sites available tends to decrease through time: there is an order of magnitude more pollen data available for the middle Holocene (ca 6000 yr BP) than for the LGM, for example \citet{harrison2016what}.
Ideally, a benchmark data set for model evaluation would provide continuous climate fields.
However, while gridding the data sets at a scale comparable to that of the climate models (see e.g. \citealp{bartlein2011pollen}) can improve the situation, this still does not solve the problem of significant gaps in site-based data coverage.

Alternative approaches to generating spatially continuous palaeoclimate reconstructions have been developed that involve combining observations with model simulations of palaeoclimates.
\citet{goosse2006using}, for example, used observations to select the most realistic member from an ensemble of climate-model simulations.
They ran a relatively large ensemble of simulations using a range of different climate forcings and/or model parametrisations to encompass uncertainties, and then selected the members of the ensemble that best matched the observations at each time step before running these simulations for longer to gain an new estimate of the climate.
In this approach, the most realistic climate is taken to be the simulated climate(s) that best matched observations after multiple simulations.
Although this approach provides continuous and self-consistent fields of climate variables, the reconstructions cannot deviate fundamentally from the model predictions and thus could still be influenced by systematic errors inherent in the model construction.
\citet{annan2013new} also used an ensemble of model simulations, but in this case they used multiple models.
The ultimate climate reconstruction was assumed to be a weighted average of those climate models, where the weighting was determined by the goodness-of-fit to the observations.
They applied a global weighting to each model rather than allowing the goodness-of-fit to vary regionally.
As a result, there are regions where the reconstructed palaeoclimate is far from the observations, producing a palaeoclimate reanalysis that is highly influenced by systematic errors in the models.

Variational data assimilation techniques provide a way of combining observations and model outputs to produce climate reconstructions that are not explicitly constrained to a given source \citep{nichols2010mathematical, lahoz2014data}.
Variational techniques are widely used by the weather forecasting community (e.g. \citealp{daley1994atmospheric}) and have also been used to reconstruct palaeoclimate. 
\citet{gebhardt2008reconstruction}, for example, used this approach to reconstruct European climates during the Last Interglacial. 
\citet{simonis2012reconstruction} applied the same basic approach to reconstruct January and July temperatures across European climate during the late Glacial (13,000 yr B.P.) and early Holocene (8,000 yr B.P.). 
The method involves applying a spatial constraint, based on a two-dimensional advection-diffusion equation of atmospheric dynamics, to upscale climate variables derived from statistical transfer functions relating the abundance of plant taxa with January and July temperature. 
In both examples, modern-day wind fields were used as the prior to determine the spatial scale and assumed to be the same in the past.

\citet{tardif2018last} also use variational techniques to create palaeoclimate reconstructions for the Last Millennium, using an ensemble of transient palaeoclimate simulations.
They first determine the relationship between palaeoclimate reconstructions and the model-derived prior using linear regression, and then determine the best linear unbiased estimate (BLUE) using the Kalman formulation, to create the analytical reconstructions. 
Thus, temporal relationships are not based on an explicit analytical function designed to preserve structures (auto-correlations and/or discontinuities) in the prior and/or observations.  
Spatial correlation is generated from the prior ensemble with a covariance localisation applied to prevent spurious correlations. 
This, and the necessity to define scaling parameters, involves a number of arbitrary choices which influence the final reconstructions and make it difficult for these reconstructions to deviate substantially from the prior.

The 3D-Variational method finds the maximum a posteriori Bayesian estimate of the palaeoclimate given the site-based reconstructions and a prior estimate.
While this could lead to the generation of reconstructions with sharp changes in time and/or space, it is possible to incorporate additional assumptions about the error of the prior estimate (the difference between the true climate and the prior) to prevent this by ensuring continuity of the solution.
The degree of continuity in the change of the reconstructed climate field can be controlled by adjusting two length scales: a spatial length scale that determines how smooth the spatial correlation in the prior is between different geographical areas and a temporal length scale that determines how smooth it is through the seasonal cycle.

Here we apply this method to reconstruct six palaeoclimate variables across southern Europe at the LGM.
The six climate variables are those provided in the \citet{bartlein2011pollen} dataset, namely mean annual temperature (MAT, $^{\circ}C$), mean temperature of the coldest month (MTCO, $^{\circ}C$), mean temperature of the warmest month (MTWA, $^{\circ}C$), growing degree days above a baseline of above $5^{\circ}$C (GDD5, $d ^{\circ}C$), mean annual precipitation (MAP, mm) and an index of plant-available moisture (the ratio of actual to equilibrium evapotranspiration or $\alpha$ in \citet{bartlein2011pollen} re-expressed as a moisture index (MI, unitless) defined as the ratio of MAP to equilibrium evapotranspiration in our analyses. 
The conversion was made using the \citet{zhang2004rational} formulation of the Budyko relationship).
We use pollen-based reconstructions of climatic variables for the region of southern Europe (defined here as south of 50$^{\circ}$N and extending eastward to 50$^{\circ}$E) from \citet{bartlein2011pollen} as our observations.  
Although the sites from Europe were used to produce a gridded map in \citet{bartlein2011pollen}, here we used the underlying individual site reconstructions. 
Some of the reconstructions used in \mbox{\citet{bartlein2011pollen}} were derived by model inversion, and these were excluded from our data set. 
\citet{bartlein2011pollen} gives mean values as anomalies from the modern climate, as well as standard errors.
We use eight LGM climate simulations (CCSM4, CNRM-CM5, MPI-ESM-P, MRI-CGCM3, FGOALS-g2, COSMOS-ASO, IPSL-CM5A-LR, MIROC-ESM) from the 3rd phase of the Palaeoclimate Modelling Intercomparison Project (PMIP3: \citealp{braconnot2012evaluation}) to create a prior.
These simulations were forced by changes in incoming solar radiation, changes in land-sea geography and the size and extent of ice sheets, and a reduction in atmospheric [CO$_2$] (see \citealp{braconnot2012evaluation} for details of the modelling protocol).

Our approach introduces features novel to the field of palaeoclimate data assimilation, explicitly designed to maximise the usefulness of the reconstructions for climate model evaluation. 
Specifically, by solving the full variational problem we take into account nonlinearities in the system.
Furthermore we minimise the dependency of the final analytical reconstructions on the prior generated from the climate models by using a prescribed correlation function for the error of the prior and by using a resolution matrix \citep{menke2012geophysical, delahaies2017constraining} to determine the temporal correlation length scale. 
The resolution matrix provides a particularly useful way to overcome problems caused by the sparsity of site-based palaeoclimate reconstructions at the LGM. 
In addition to investigating methods to determine appropriate spatial and temporal length scales, we provide a way of calculating the error in the final reconstructions.

\section{Data Assimilation with Spatial and Temporal Correlations in the Prior}
\label{sec:data_assimilation}
In this section we describe the underlying method used in this paper.
Section \ref{sec:inverse_problem} describes the inverse problem solved by the method and the types of data used.
Section \ref{sec:background_error_correlation} shows how the different climate variables can be related to one another by specifying correlations from our prior estimate of the system. 
Finally section \ref{sec:preconditioning} describes how the problem is preconditioned in order to reduce the computation cost.

\subsection{The Inverse Problem}
\label{sec:inverse_problem}
Our problem is to determine the palaeoclimate that existed from a particular set of site-based reconstructions.
We label the reconstructions as the column vector $\bo{y}_i \in \mathbb{R}^6$ for site $i$.
For each reconstruction, $\bo{y}_i$, there are a total of $6$ variables that may have been reconstructed, namely; $\alpha$, MAP, MAT, MTCO, MTWA and GDD5. 
All these reconstructions together make the observations labelled $\bo{y} \in \mathbb{R}^{6N}$ such that
\begin{linenomath*}
    \begin{equation}
        \bo{y} = \left(\bo{y}_1^T | \bo{y}_2^T | \cdots | \bo{y}_N^T \right)^T
    \end{equation}
\end{linenomath*}
where $N$ is the number of reconstructions.
The reconstruction technique gives the variances for each reconstruction that we label as the column vector $\bo{v}_y \in \mathbb{R}^{6N}$ in the same order as $\bo{y}$.
Not all variables are reconstructed at every site, for these variables the variance tends to infinity; this is achieved by setting their inverse to 0. 

From these reconstructions we want to produce a gridded climate, the state vector, $\bo{x} \in \mathbb{R}^{13M}$ where there are $M$ grid cells.
The $j$'th grid cell of the state is labelled $\bo{x}_j \in \mathbb{R}^{13}$ where 
\begin{linenomath*}
    \begin{equation}
        \bo{x} = \left(\bo{x}_1^T | \bo{x}_2^T | \cdots | \bo{x}_M^T \right)^T.
    \end{equation}
\end{linenomath*}
For each grid cell the method determines a set of $13$ variables: the mean annual precipitation ($P$) and the $12$ average temperatures for each month, $\bo{T}$ where 
\begin{linenomath*}
    $$ \bo{T} = \left( T_1 \: T_2 \: \dots \: T_{12}\right)^T$$
\end{linenomath*}
where $T_m$ is the temperature at month $m$. 

For a general function $\bo{h}$ that maps a gridded climate $\bo{x}$ to the site-based observations we state the problem as trying to find an $\bo{x}$ such that
\begin{linenomath*}
    \begin{equation}
        \label{eq:uncorrected_obs_function}
        \bo{h} ( \bo{x} ) = \bo{y}.
    \end{equation}
\end{linenomath*}
Solving equation \eqref{eq:uncorrected_obs_function} for $\bo{x}$ is ill-posed as there are several $\bo{x}$ that are possible solutions. 
A prior estimate of the state called the background or prior ($\bo{x}_b$) allows us to find the best $\bo{x}$ that solves equation \eqref{eq:uncorrected_obs_function} and remains close to the prior.
The standard deviations of the prior are labelled as the vector $\bo{v}_b \in \mathbb{R}^{13M} $ in the same order as $\bo{x}_b$.

It can be shown \citep{nichols2010mathematical} that the optimal solution of equation \eqref{eq:uncorrected_obs_function} with a prior estimate of the state is defined as the analysis, $\bo{x}_a$, where
\begin{linenomath*}
    \begin{equation}
        \label{eq:analysis_def}
        \bo{x}_a = \min_{\bo{x}} J(\bo{x}).
    \end{equation}
\end{linenomath*}
with the cost function $J$ as 
\begin{linenomath*}
    \begin{equation}
        \label{eq:dimensional_cost_function}
        J(\bo{x}) = \frac{1}{2} (\bo{x} - \bo{x}_b)^T \bo{B}^{-1} (\bo{x}-\bo{x}_b) + \frac{1}{2}(\bo{y} - \bo{h}(\bo{x}))^T \bo{R}^{-1} (\bo{y} - \bo{h}(\bo{x})).
    \end{equation}
\end{linenomath*}
Here $\bo{B}$ is the covariance of the uncertainties in the prior (conventionally denoted $\bo{B}$, for background) and $\bo{R}$ is the covariance of the uncertainties in the site-based reconstructions. 
Equations \eqref{eq:analysis_def} and \eqref{eq:dimensional_cost_function} ensure that the solution is the optimal distance from the observations subject to ensuring that the solution is not too far from the prior estimate, weighted by the error statistics in each.
We assume that there are no correlations in the errors of the observations so we set
\begin{linenomath*}
    $$\bo{R} = diag(\bo{v}_y).$$
\end{linenomath*}
The prior error covariance matrix can be represented as the product of the standard deviations of the prior and the correlations between the errors in the variables in the prior. 
Hence we write 
\begin{linenomath*}
    \begin{equation}
        \label{eq:B_covariance}
        \bo{B} = \bo{\Sigma} \bo{C} \bo{\Sigma}
    \end{equation}
\end{linenomath*}
where 
\begin{linenomath*}
    \begin{equation}
        \bo{\Sigma} = diag(\bo{v}_b^{\frac{1}{2}}),
    \end{equation}
\end{linenomath*}
is the diagonal matrix formed of the standard deviations of the prior error and $\bo{C}$ is the prior error correlation matrix.

\subsection{Prior Error Correlation}
\label{sec:background_error_correlation}
The difference between the true $\bo{x}$ and the prior, the error in the prior, is expected to be smooth between adjacent grid cells and also from month to month since it would be unlikely that the observations would contain sharp jumps in climate that aren't present in the prior. 
It would be unusual, for example, to have very high temperature in March if the temperatures in February and April are very low, if this behaviour isn't seen in the prior.
To achieve this we impose a structure on the prior error correlation matrix, $\bo{C}$, that weighs the cost function so that its minimum is smooth. 
This allows the prior error to be smooth, but still allows non-smooth areas if there is significant evidence to support it in the prior and/or the observations.

We assume there are two independent sets of correlations in the prior.
The first correlation is spatially between the different grid cells.
We also assume that the spatial correlation between the grid cells is homogeneous and valid on a sphere, so that for an angle $\theta_{ij}$ on a great circle of the Earth between each cell $i$ and $j$ the correlation is given by,
\begin{linenomath*}
    \begin{equation}
        \label{eq:gaussian_correlation_function}
        c_{L}(\theta_{ij}) = \left( \frac{a}{L} \sin\left(\frac{\theta_{ij}}{2}\right) \right) \mathcal{K}_1\left( \frac{a}{L} \sin\left(\frac{\theta_{ij}}{2}\right) \right)
    \end{equation}
\end{linenomath*}
where \mbox{$c_L$} is a case of a Mat\'ern function \citep{matern1986spatial, handcock1994approach} with order 1 and $\mathcal{K}$ is the modified Bessel function of the second kind, evaluated using the boost C++ library \citep{boost2018maddock}.
Here the correlation length scale is $L=L_s$ and $a=6371$km is the radius of the Earth.
The correlation matrix between all grid cells, $\bo{C}_{L_s}$, is given as 
\begin{linenomath*}
    $$(\bo{C}_{L_s})_{ij} = c_{L_s}(\theta_{ij}).$$
\end{linenomath*}
The choice of $L_s$ is dependent on the datasets used in $\bo{y}$ and $\bo{x}_b$ and so is specific to each problem. 
In section \ref{sec:finding_length_scales} a method of finding $L_s$ is shown for a particular experiment. 

The second assumed correlation is between the error in the average temperatures of the prior.
We assume that there is a correlation between the average temperatures of a month and the surrounding months given by equation \eqref{eq:gaussian_correlation_function}. 
Here $\theta_{ij} = mod_{12}(\left| i - j \right|)$ between months $i$ and $j$.
The correlation length scale is $L = L_t$ and $a=6 / \pi$.
The appropriate value of $L_t$ again depends on the datasets given and is shown for a particular experiment in section \ref{sec:finding_length_scales}.
For each grid cell the correlation between the different climate variables is given by $\bo{C}_{L_t}$ where 
\begin{linenomath*}
    \begin{equation} 
        \label{eq:temp_cor_mat}
        \bo{C}_{L_t} = \left(
        \begin{array}{c | ccc}
            1      & 0 & \dots                 & 0 \\
            \hline
            0      &   &                       &   \\
            \vdots &   & \{c_{L_t}(\theta_{ij}\}_{ij} &   \\
            0      &   &                       &   \\
        \end{array}
        \right).
    \end{equation} 
\end{linenomath*}
Note how $\{c_{L_t}(\theta_{ij}) \}_{ij}$ is offset by the first row and column due to the presence of the precipitation term which is uncorrelated to the temperature terms. 

These two sets of correlations imply that all the variables in the error of the prior are correlated.
For instance the grid cells $i$ and $j$ are correlated by $\left( \bo{C}_{L_s} \right)_{ij}$ and the temperatures in month $l$ and $k$ are correlated by $\left(\bo{C}_{L_t} \right)_{lk}$.
This means that the temperatures in month $l$ in grid cell $i$ and month $k$ in grid cell $j$ are correlated by the product $\left( \bo{C}_{L_s} \right)_{ij} \left( \bo{C}_{L_t} \right)_{lk}$. 
Repeating this for every variable gives an overall correlation for the prior ($\bo{C}$ from equation \eqref{eq:B_covariance}) as
\begin{linenomath*}
    \begin{equation}
        \bo{C} = \bo{C}_{L_s} \otimes \bo{C}_{L_t}
    \end{equation}
\end{linenomath*}
where $\otimes$ is the Kronecker product of matrices.

The incorporation of correlations structures is due to the fact that the state space covers space and time.
We introduce the $\bo{C}_{L_s}$ and $\bo{C}_{L_t}$ correlations to make the prior error smooth in space and time respectively. 
The presence of the scales $L_s$ and $L_t$ allows the adjustment of the smoothing in both dimensions and should depend, at least in part, on the spatial and temporal distribution of the prior and site-based reconstructions. 
In section \mbox{\ref{sec:finding_length_scales}} we discuss methods for choosing these scales.

\subsection{Preconditioning and the Condition Number}
\label{sec:preconditioning}
The minimum of the cost function is sensitive to change in the input data of the problem and to computational errors. 
This sensitivity reflects the difficulty in solving the problem and is measured by the condition number of the Hessian of the cost function \citep{golub1996matrix}. 
We define the condition number $\kappa$ of a symmetric positive definite matrix $\bo{M}$ to be 
\begin{linenomath*}
    \begin{equation}
        \label{eq:condition_number_definition}
        \kappa(\bo{M}) = \frac{\lambda_{\text{max}(\bo{M})}}{\lambda_{\text{min}(\bo{M})}}
    \end{equation}
\end{linenomath*}
where $\lambda_{\text{max}(\bo{M})}$ and $\lambda_{\text{min}(\bo{M})}$ are the maximum and minimum eigenvalues of $\bo{M}$.
Here, $\bo{M}$ is the Hessian of the cost function, given by its (first order) second derivative $\bo{S} = \bo{H}\bo{B}\bo{H}^T + \bo{R}$.
This condition number indicates the computational work needed to minimise the cost function.
Equation \eqref{eq:condition_number_definition} shows how the condition number of $\bo{S}$ represents the disparity in scales of the problem. 
As the eigenvalues represent the sizes of the scales of $\bo{S}$, their ratio represents the largest scale that will be encountered when inverting $\bo{S}$. 
Since large scale differences create more numerical inaccuracy, a large condition number will increase the computational cost and lead to an inaccurate solution.

\citet{haben2010conditioning} shows that the bounds on the condition number can be reduced by minimising the cost function around $\bo{w}$ instead of $\bo{x}$ where
\begin{linenomath*}
    \begin{equation}
        \label{eq:x_to_w_transform}
        \bo{B}^{\frac{1}{2}} \bo{w} = \bo{x} - \bo{x}_b
    \end{equation}
\end{linenomath*}
where $\bo{B}^{\frac{1}{2}}$ is the symmetric square root of the matrix $\bo{B}$ such that 
\begin{linenomath*}
    $$ \bo{B} = \bo{B}^{\frac{1}{2}} \bo{B}^{\frac{1}{2}}.$$
\end{linenomath*}
The use of this linear transformation can be thought of as a z-score to work with uncorrelated states.

Equation \eqref{eq:x_to_w_transform} transforms the inverse problem from equation \eqref{eq:analysis_def} into finding
\begin{linenomath*}
    \begin{equation}
        \bo{w}_a = min_{\bo{w}} J(\bo{w}).
    \end{equation}
\end{linenomath*}
where $J(\bo{w})$ is
\begin{linenomath*}
    \begin{equation}
        \label{eq:translated_cost}
        J(\bo{w}) = 
        \frac{1}{2} \bo{w}^T \bo{w} + \frac{1}{2} (\bo{y} - \bo{h}(\bo{x}_b + \bo{B}^{\frac{1}{2}} \bo{w}) )^T \bo{R}^{-1} (\bo{y} - \bo{h}(\bo{x}_b + \bo{B}^{\frac{1}{2}} \bo{w})).
    \end{equation}
\end{linenomath*}
We use the limited memory Broyden-Fletcher-Goldfarb-Shanno (L-BFGS) method to find the state, $\bo{w_a}$, which has the minimum $J$, L-BFGS is a quasi-Newton method that maintains a limited memory version of an approximated Hessian as described in \citet{liu1989limited}. 
At each evaluation step we calculate the gradient of $J$ as
\begin{linenomath*}
    \begin{equation}
        \label{eq:conditioned_cost_jacobian}
        \nabla J(\bo{w}) = \bo{w} - \bo{B}^{\frac{1}{2}} \bo{H}^T_{\bo{x}} \bo{R}^{-1} \left(\bo{y} - \bo{h}_u(\bo{x}_b + \bo{B}^{\frac{1}{2}} \bo{w}) \right) 
    \end{equation}
\end{linenomath*}
where $\bo{H}_{\bo{x}}$ is the Jacobian of $\bo{h}$ evaluated at $\bo{x}$.
Once $\bo{w}_a$ is found we use equation \eqref{eq:x_to_w_transform} to find $\bo{x}_a$, the solution.

The error covariance of the analysis, $\bo{x}_a$, is denoted by $\bo{A}$ and is given (to first order) following \mbox{\citet{nichols2010mathematical}} as
\begin{linenomath*}
    \begin{equation}
        \bo{A} = \left( \bo{I} - \bo{K} \bo{H}_{\bo{x}_b} \right) \bo{B}.
    \end{equation}
\end{linenomath*}
where the gain matrix $\bo{K}$ is 
\begin{linenomath*}
    \begin{equation}
        \label{eq:gain_matrix}
        \bo{K} =  \bo{B} \bo{H}^T_{\bo{x}_b} \left( \bo{H}_{\bo{x}_b} \bo{B} \bo{H}^T_{\bo{x}_b} + \bo{R} \right)^{-1}.
    \end{equation}
\end{linenomath*}

\section{Experimental Design}
\label{sec:example}
We use our method to reconstruct the palaeoclimate of southern Europe during the Last Glacial Maximum (LGM).
The LGM had insolation forcing relatively similar to the present day but northern hemisphere ice sheets were more extensive, sea-level was lower and the area of the continents therefore larger, and the atmospheric [CO$_2$] was less than half of the concentration today. 
In this section we describe the choices of $\bo{h}$, $\bo{y}$ and $\bo{x}_b$ used to make this reconstruction and our choices for $L_t$ and $L_s$, the correlation length scales.

\subsection{Experiment Setup}
We use pollen-based reconstructions of climatic variables from \citet{bartlein2011pollen} as our observations. 
\citet{bartlein2011pollen} gives means as anomalies from the modern climate as well as standard errors.
We add the anomalies to the CRU CL v2.0 dataset \citep{new2002high} to derive absolute climate reconstructions.
We non-dimensionalise the climate variables in order to avoid computational issues because they are on different scales in the calculation of the cost function. 
After solving for the non-dimensionalised case, we re-dimensionalise each of the variables to be on the original scale of the observations and the prior. 
Details of the dimensionalisation and non-dimensionalisation of variables can be found in Appendix \ref{sec:non-dimensionalisation}.
We use the non-dimensionalised variables as our $\bo{y}$ and their non-dimensional standard errors, formed from the product of the standard errors and the derivative of $D_y$ (equation \ref{eq:obs_non_dim}), as $\bo{v}^{\frac{1}{2}}_{\bo{y}}$.

We use the LGM outputs from PMIP3 as our prior. 
We use the variables of monthly precipitation (that are summed to annual precipitation), monthly temperature and monthly total cloud fraction.
For each of the selected PMIP models that ran an LGM experiment we interpolate the output to a $2^{\circ} \times 2^{\circ}$ grid producing a set of maps all at the same resolution.
In order to minimise the impact of potential individual systematic model biases the simulated climate at the LGM, experiments are generally expressed relative to that specific model's pre-industrial control (PI) experiment.
We therefore interpolate each of the PI experiments to the same grid and take the difference between the LGM and PI experiments of each model as the anomaly to the modern day for each model. 
We then sum each model's anomalous values with values from the modern day (from CRU CL v2.0, as above, bilinearly interpolated to the $2^{\circ} \times 2^{\circ}$ grid) in order to produce absolute values for each model.
For each variable in the set we take the mean and variance across the set of all models to produce a gridded map.
As for the observation space, we non-dimensionalise the state space to remove any dimensional effects using equation \eqref{eq:state_non_dim}. The non-dimensional variables form the prior $\bo{x}_b$ and their non-dimensional variances, formed from taking the product of the variances and the derivative of equation \eqref{eq:state_non_dim}, form $\bo{v}_{\bo{x}_b}$.

The observation function, $\bo{h}$, links together the variables from both datasets. 
At each site, $i$, we define the observation function as 
\begin{linenomath*}
    \begin{equation}
        \label{eq:observation_function_geo}
        \bo{\hat h} (\bo{x}_i) =
        \bo{\hat h}\left(
        \begin{array}{c}
            \bar{P} \\
            \bar{\bo{T}} \\
        \end{array}
        \right) =
        \left(
        \begin{array}{c}
            \mu(\bo{x}_i) \\
            \bar{P}\\ 
            mean(\bar{\bo{T}})\\ 
            max(\bar{\bo{T}})\\
            min(\bar{\bo{T}})\\ 
            G(\bar{\bo{T}})
        \end{array}
        \right).
    \end{equation}
\end{linenomath*}
The derivatives, $\frac{\partial max(T)}{\partial T_m}$ and $\frac{\partial max(T)}{\partial T_m}$ are taken to be $1$ if $T_m$ is the maximum or minimum of $T$ and $0$ elsewhere.
The moisture index function $\mu$ is 
\begin{linenomath*}
    \begin{equation}
        \label{eq:alpha}
        \mu(\bo{x}_i) = 1 + m(\bo{x}_i) - \left( 1 + m(\bo{x}_i)^{\omega} \right)^{\frac{1}{\omega}}
    \end{equation}
\end{linenomath*}
as given by the Budyko curve with $\omega =3$ as described in \citet{zhang2004rational}.
The moisture index $m$ is calculated as
\begin{linenomath*}
    \begin{equation}
        \label{eq:mi_true_calc}
        m(\bo{x}_i) = P \lambda \left[ \sum_k^{12} {l}_k \frac{R(T_k, S_k) \frac{\partial e_s}{\partial T}\big|_{T_k}}{\frac{\partial e_s}{\partial T}\big|_{T_k} + \gamma} \right]^{-1}
    \end{equation}
\end{linenomath*}
where $\gamma$ ($0.067 kPaK^{-1}$) is the psychrometer constant at sea level, $l_j$ is the length of month $j$ in days and where  
\begin{linenomath*}
    $$ \frac{\partial e_s}{\partial T} = \frac{10.5485}{\left(237.3 + T\right)^2} exp \left(\frac{17.27 T}{237.3 + T} \right), $$
\end{linenomath*}
is the differentiated Roche-Magnus formula from \citet{allen1998crop}.
The function $R(T_k, S_k)$ is the daily net radiation at the vegetated surface defined in \citet{davis2016simple} for the middle day in month $k$. 
The variable $S_k$ is the total cloud fraction for month $j$ which is taken from the PMIP3 average described above.
We define 
\begin{linenomath*}
    $$
    G(\bar{\bo{T}}) = \frac{1}{N_y} \sum_k^{12}
    \begin{cases}
        l_k \left( \bar{T}_k - \frac{5}{T_s} \right) & \bar{T}_k > \frac{5}{T_s} \\
        0 & \text{else}
    \end{cases},
    $$
\end{linenomath*}
and the mean function to be $mean(\bar{\bo{T}}) = \frac{1}{N_y} \sum_k^{12} l_k \bar{T}_k$ and $max(\bar{\bo{T}})$ and $min(\bar{\bo{T}})$ to be the maximum and minimum temperature in $\bar{\bo{T}}$ respectively.
The full observation function, $\bo{h}$, is formed by applying $\hat{\bo{h}}$ at each grid cell where there is an observation and defining 
\begin{linenomath*}
    $$\bo{h}(\bo{x}) = \left( \hat{\bo{h}}(\bo{x}_1)^T | \hat{\bo{h}}(\bo{x}_2)^T | \cdots \right) $$
\end{linenomath*}
and so $\bo{h}$ will have the dimension $6N$, and hence the Jacobian of $h$, $H$, will have dimension $(6N)^2$.

\subsection{Determining $L_t$ and $L_s$}
\label{sec:finding_length_scales}
The two correlation length scales, $L_t$ and $L_s$, in $\bo{C}$ (section \ref{sec:background_error_correlation}) determine the strength of the correlation in the errors of the prior.
By varying the length scales we can vary how smooth the error of the prior is and hence how smooth the solution is. 
If the length scale is too large then the error will be over-smoothed and the solution will miss smaller scale features such as inter-annual temperature changes or spatially small features such as topography.
A length scale too small will mean the solution will be too coarse and contain unrealistic jumps.

\begin{figure*}[h]
    \centering
    \includegraphics[width=\textwidth]{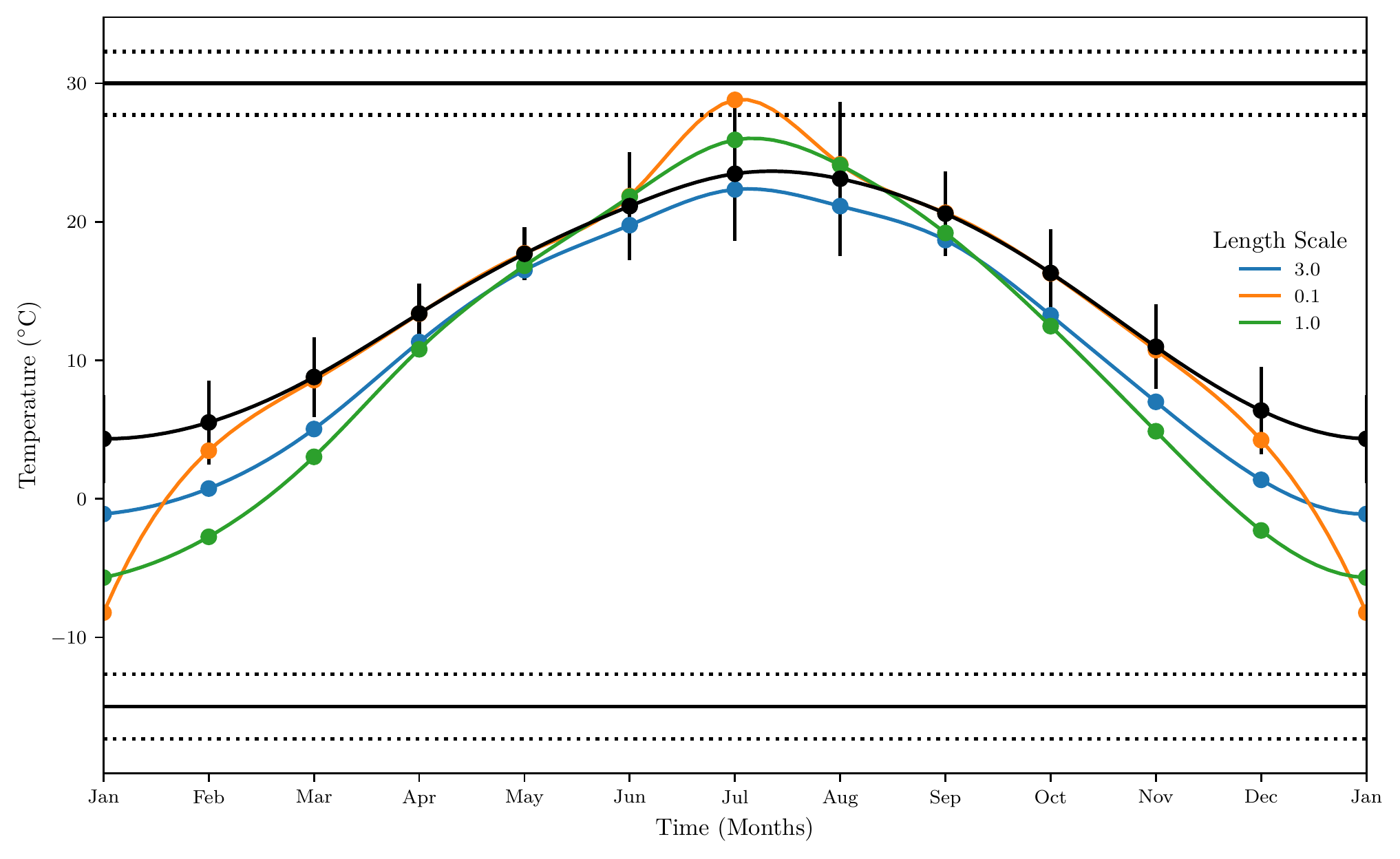}
    \caption{
        \label{fig:L_t_test}
        Yearly temperature for the assimilation performed on a single simulated site at N$37.50^{\circ}$ and E$33.73^{\circ}$ with varying values of $L_t$. 
        The different coloured dots are the results of the assimilation for different values of $L_t$.
        The black dots in the centre are the prior for the grid cell that contains this site with error bars of 1 standard deviation. 
        The B-spline interpolation of the dots is shown as the curved lines.
        The observations of MTWA and MTCO are represented by the higher and lower solid black lines respectively with the dotted lines showing 1 standard deviation around the their mean.}
\end{figure*}

In order to determine a suitable value for $L_t$ we consider a single grid cell with a single simulated observation at $37.50^{\circ}$N and E$33.73^{\circ}$, which allows us to ignore the effects of $\bo{C}_{L_s}$.
The example only has observations of MTCO and MTWA ($-15^{\circ}$C and $30^{\circ}$C respectively), allowing us to ignore the non-linear effects of calculating $\alpha$.
Fig. \ref{fig:L_t_test} shows the prior and observations for the sample as well as the estimated states after assimilation for different values of $L_t$.
For all values of $L_t$ the analysis doesn't match the observed MTCO since the prior temperature for January has low uncertainty.
Low values of $L_t$ create an analysis that swaps between the prior and the observations.
Although the solution always matches either the reconstructions or the prior, the jumps between them are unrealistic.
On the other hand high values of $L_t$ create an analysis that follows the prior too closely and is unable to create high and low temperatures.
The value of $L_t = 1$ produces an assimilation that follows the shape of the prior but lies between the values of the prior and the observations.

We can further understand $L_t$ by seeing how information is changed by the method.
If we consider the hypothetical, true solution to the inverse problem, $\bo{w}_t$, then by equation \eqref{eq:uncorrected_obs_function} we have that 
\begin{linenomath*}
    $$ \bo{H}_{\bo{x}_b} \bo{B}^{\frac{1}{2}} \bo{w}_t \approx \bo{y} - \bo{h}(\bo{x}_b)$$
\end{linenomath*}
since, up to first order,
\begin{linenomath*}
    \begin{equation}
        \bo{H}_{\bo{x}_b} (\bo{x} - \bo{x}_b) \approx \bo{h}(\bo{x} - \bo{x}_b).
    \end{equation}
\end{linenomath*}
Further \citet{nichols2010mathematical} shows how 
\begin{linenomath*}
    $$\bo{x}_a - \bo{x}_b \approx \bo{K} \left( \bo{y} - \bo{h} (\bo{x}_b) \right),$$
\end{linenomath*}
 where $\bo{K}$ is the gain matrix defined in equation \eqref{eq:gain_matrix}.
Hence we can consider the change from true solution to our computed one ($\bo{w}_a$) as being given by 
\begin{linenomath*}
    $$\bo{w}_a \approx \bo{N} \bo{w}_t$$
\end{linenomath*}
where 
\begin{linenomath*}
    $$\bo{N} = \bo{B}^{-\frac{1}{2}} \bo{K}  \bo{H}_{\bo{x}_b} \bo{B}^{\frac{1}{2}} $$
\end{linenomath*}
is the resolution matrix as described in \citet{menke2012geophysical, delahaies2017constraining}.

Resolution matrices where the diagonal elements are close to 0 describe a situation where, if perfect information is input, then the solution would only contain part of this information.
In situations where the resolution matrix has large off-diagonal terms, the solution is degraded by interference between variables.
If the opposite is true, the resolution matrix is close to the identity matrix.
The best method will have a resolution matrix that resolves as many variables as possible whilst having few variables interfering with each other. 

\begin{figure*}[h]
    \centering
    \includegraphics[width=\textwidth]{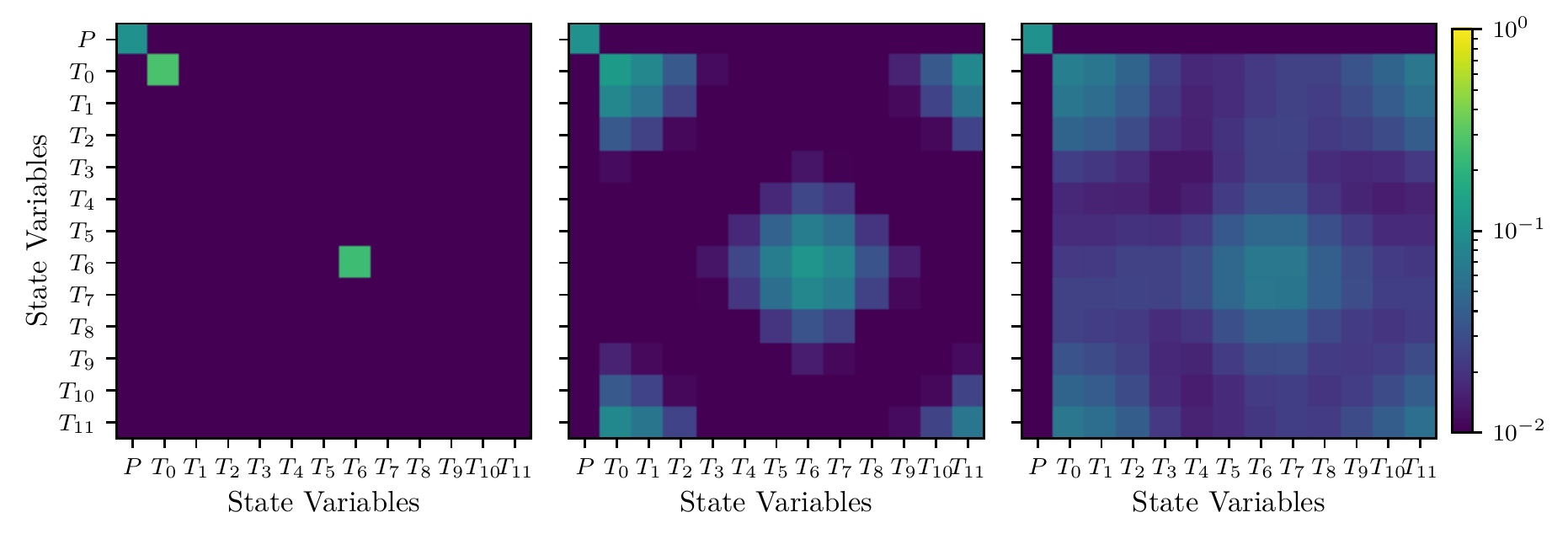}
    \caption{
        \label{fig:res_mat_comb_paper}
        The resolution matrices for the assimilation method with a sample single grid cell and a simulated observation at N$37.50^{\circ}$ and E$33.73^{\circ}$.
        The colour is the log value of the resolution matrix $N$ for values of $L_t = 0.1, 1$ and $2$ respectively.
    }
\end{figure*}

Fig. \ref{fig:res_mat_comb_paper} shows how the resolution matrix changes with respect to $L_t$ for the same test grid cell as in Fig. \ref{fig:L_t_test}.
The simulated prior temperatures are closest to the observations in January and July such that for small values of $L_t$, the method resolves temperatures in these months well.
However, for large $L_t$ the method improves the patterns away from these months whilst degrading reconstructions of January and July. 
Values of $L_t$ in between the large and small values show a mixture of both high resolution and low interference.
These results together with the results from Fig. \ref{fig:L_t_test} suggest a value of $L_t=1$ is suitable for this problem.

The choice of the other scale, $L_s$, is especially relevant for the relatively sparse dataset used here.
A higher $L_s$ represents errors in the prior being correlated even though they are far away, whereas a low $L_s$ represents errors not being highly correlated even though they are close together. 
A large $L_s$ means that information from the reconstructions could be be propagated over a large distance. 
While this is useful in maximizing the use of a geographically sparse data set, it could be unrealistic if this propagation extends too far beyond the source area for the pollen on which the site reconstructions are based (which is generally, though not always, of the order of $20 - 100$km around the site). 
In order to obtain a realistic solution whilst maximising the use of the data we choose $L_s$ such that the assumed average source area of the different sites does not overlap.

$L_s$ corresponds to the area that each observation impacts, so an increase in $L_s$ gives higher utilisation of observations.
\citet{haben2010conditioning} show that the condition number of the inverse problem is proportional to the distance between the reconstruction sites which, in this case, is proportional to $L_s$.
However, the condition number corresponds to the sensitivity of inverting the Hessian to inputs and so is inversely proportional to the computational accuracy of the problem, up to first order.
Hence, it is important to check that a choice of large $L_s$ doesn't lead to a condition number for the problem that is too large to give an accurate result.

\begin{figure*}[h]
    \centering
    \includegraphics[width=\textwidth]{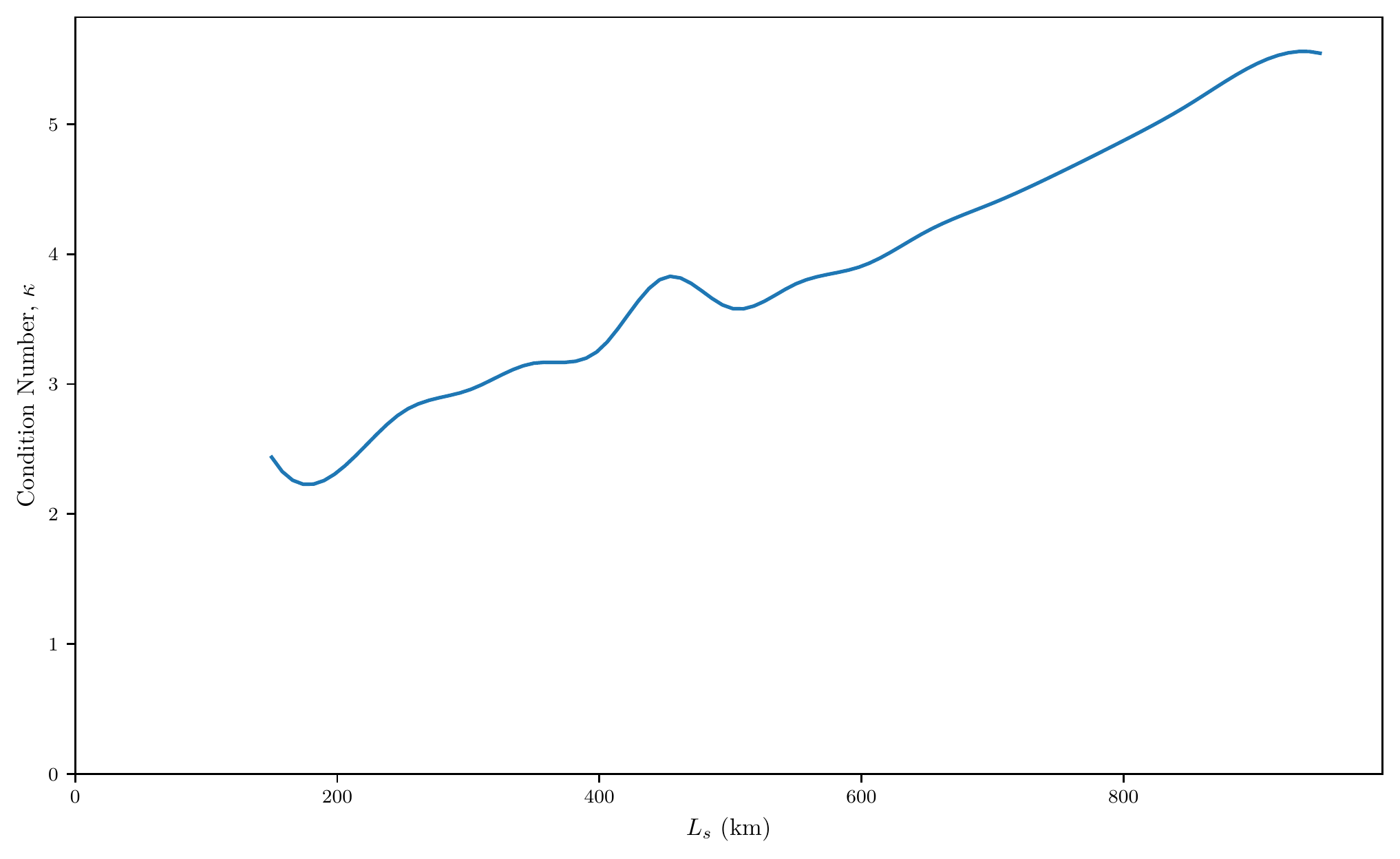}
    \caption{
        \label{fig:comb_paper_cond_num}
        The condition number of our example problem as a function of $L_s$, the spatial length scaling. 
    }
\end{figure*}

Fig. \ref{fig:comb_paper_cond_num} plots $\kappa(\bo{S})$ against $L_s$ and shows how $\kappa(\bo{S})$ begins to increase with higher $L_s$. 
Also Fig. \ref{fig:comb_paper_cond_num} shows  several inflection points which could indicate values of $L_s$ that allow multiples of observations to interact. 
For this paper we pick a value of 400km for $L_s$ as this is large enough to propagate information sufficiently far from the different reconstructions.
As seen in Fig. \ref{fig:comb_paper_cond_num}, $L_s=400$km still has a relatively low condition number and hence the solution will be relatively accurate.

\section{Results}
\label{sec:results}

\begin{figure*}[h]
    \centering
    \includegraphics[width=\textwidth]{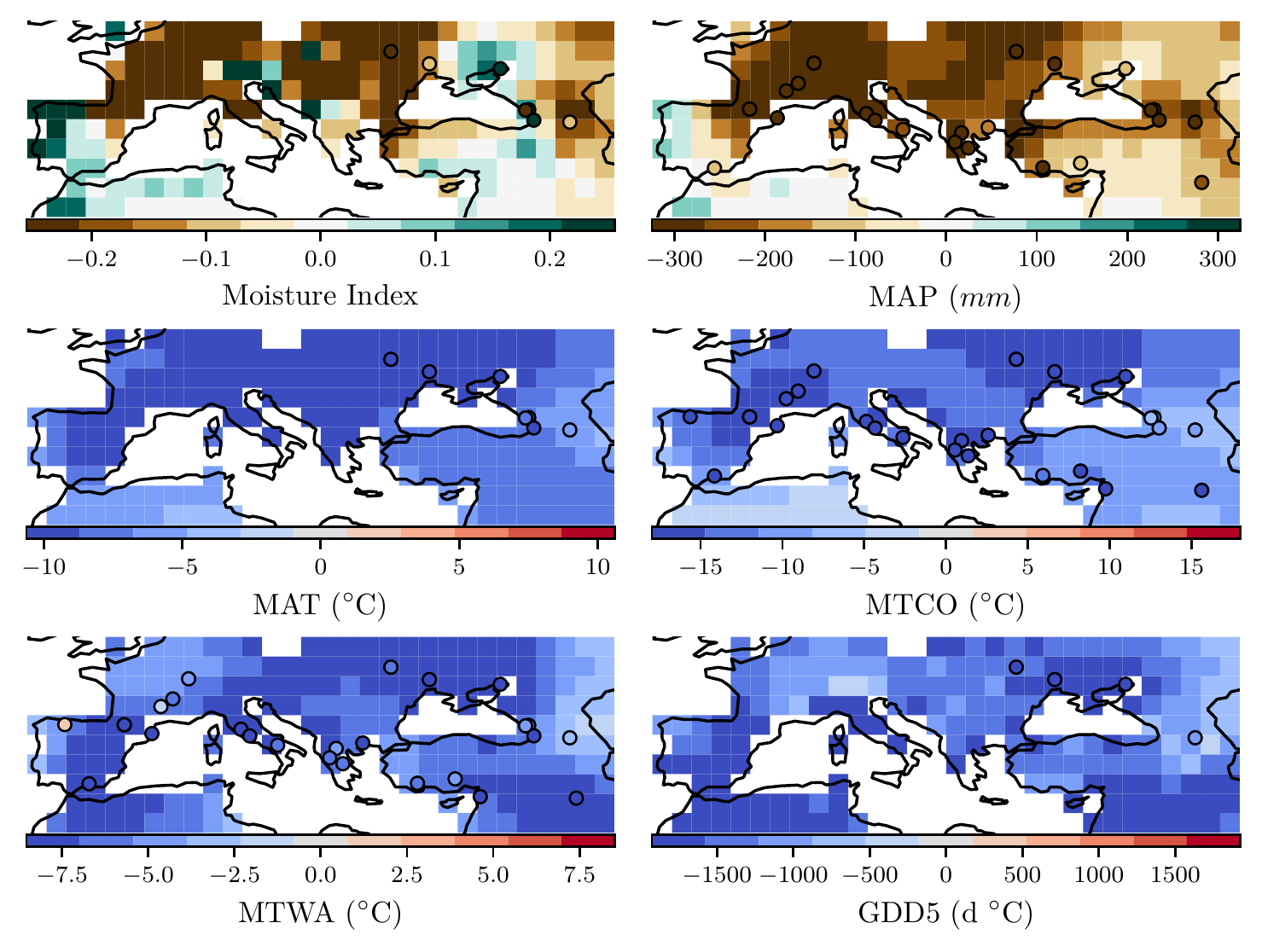}
    \caption{
        \label{fig:comb_paper_ana}
        The result, $\bo{h}(\bo{x}_a)$, is dimensionalised and represented by the colour field with the dots representing observations made ($\bo{y}$).
        Observations of $\alpha$ have been translated to moisture index through equation \eqref{eq:alpha}.
    }
\end{figure*}

The solution using scaling values of $L_t = 1$ and $L_a = 400$ (Fig. \ref{fig:comb_paper_ana}) produces climates at 50 sites and surrounding grid cells that are close to the reconstructions, as expected, over much of the region. 
However, this is not the case for the MI values of the 3 sites at the eastern tip of the Black sea (Apiancha, Kobuleti, Sukhumi).
These discrepant cases occur either where there is significant disagreement between different reconstructions and/or disagreement between the reconstructions and the prior with at least one of the reconstructions having relatively low variance. 
This reconstruction is weighted highly in the cost function and the solution does not meet the other reconstructed variables or the prior. 
This creates a situation in which the best possible solution differs from both the reconstructions and prior. 

\begin{figure*}[h]
    \centering
    \includegraphics[width=\textwidth]{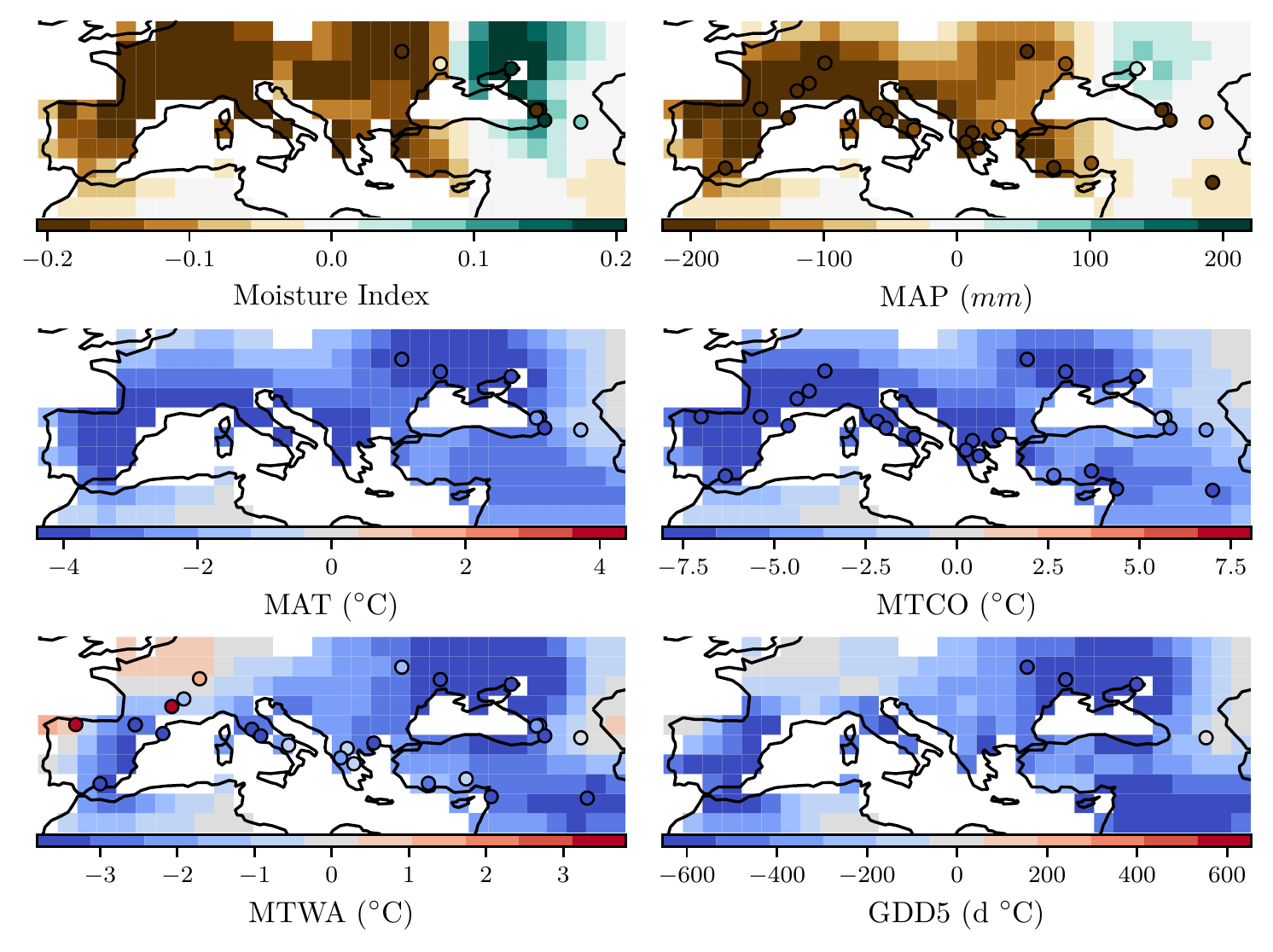}
    \caption{
        \label{fig:comb_paper_err}
        The colour field is the difference between the reconstructed climate field and the prior, $\bo{h}(\bo{x}_a) - \bo{h}(\bo{x}_b)$, dimensionalised. The dots are the differences between the site-based observations, $\bo{y}$, and the reconstructed climate of the grid cell they are in.
        Observations of $\alpha$ have been translated to moisture index through equation \eqref{eq:alpha}.
    }
\end{figure*}

The difference between the solution and the prior, transformed by equation \eqref{eq:observation_function_geo} at each grid cell and dimensionalised via equation \eqref{eq:obs_non_dim}, shows that the climate is much drier than the prior in the western part of the area, as shown by MI and precipitation (Fig. \ref{fig:comb_paper_err}).
MAT has increased in some regions but decreased in others; this suggests that the inclusion of $\bo{C}_{L_s}$ is working as intended, since although there are varied changes in MAT, the changes occur in a spatially smooth way.
Furthermore there has been an increase in temperature seasonality as MTCO has become colder at all sites and MTWA has become warmer at most sites. 
This, together with the changes to MAT and GDD5 suggests that $\bo{C}_{L_t}$ is having the desired effect; as the changes to MTCO and MTWA are impacting the whole of the seasonal cycle of the climate and giving reasonable and smooth changes to both MAT and GDD5.

\begin{figure*}[h]
    \centering
    \includegraphics[width=\textwidth]{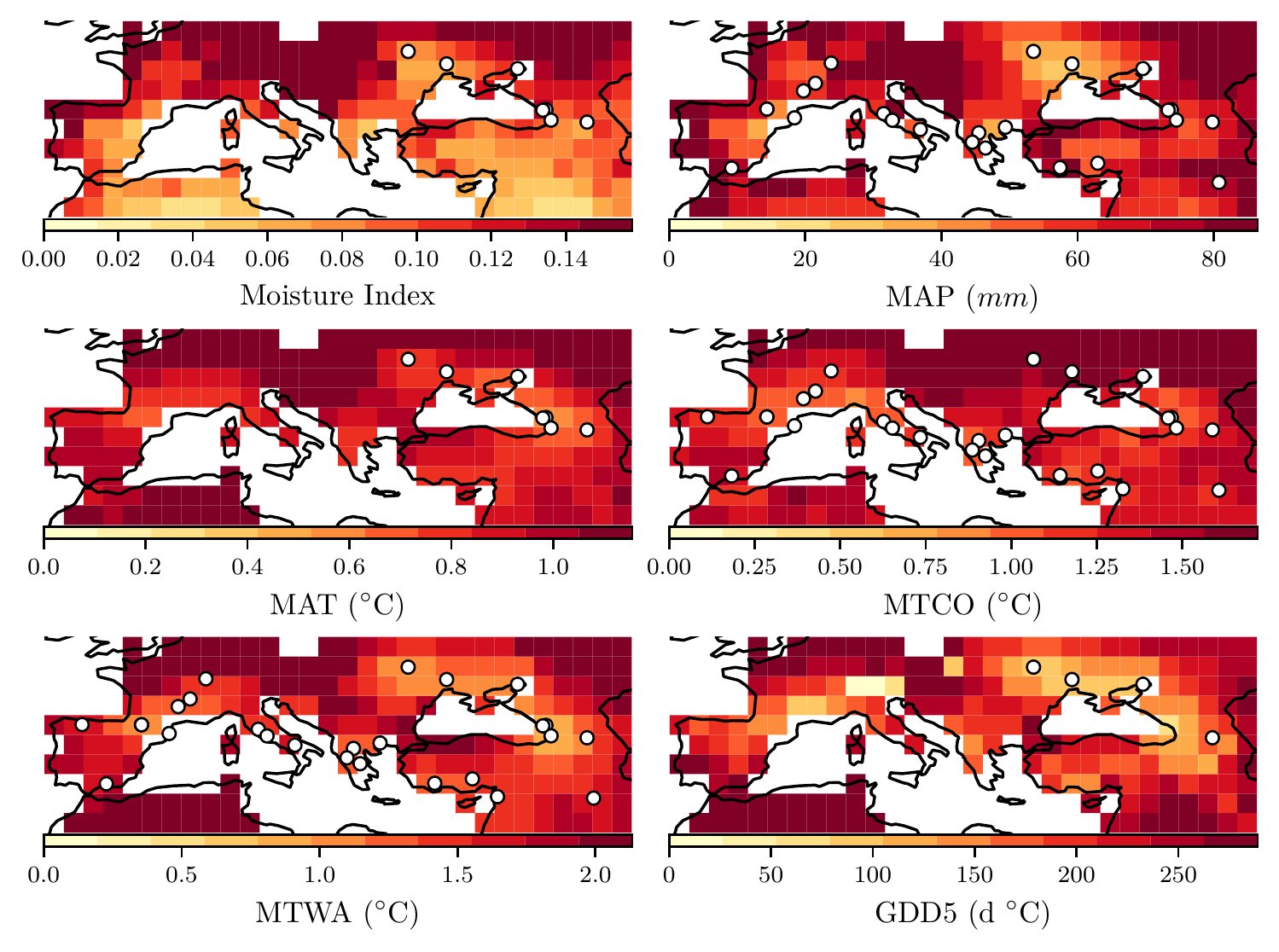}
    \caption{
        \label{fig:comb_paper_sd}
        The standard deviation of the result, given by the dimensionalised square root of the main diagonal of $\bo{H}_{\bo{x}_a} \bo{A} \bo{H}^T_{\bo{x}_a}$ (the analysis error covariance in observation space), is represented by the colour field where the dots represent sites of observations.
        Observations of $\alpha$ have been translated to moisture index through equation \eqref{eq:alpha}.
        For areas with very low temperature it is almost certain that GDD5 is zero and so these areas have been left blank.
    }
\end{figure*}

In general (Fig. \ref{fig:comb_paper_sd}) grid cells near reconstruction sites have less error, because the solution is using information from both the prior and the reconstructions, while grid cells further away from reconstruction sites have higher error by defaulting to the error in the prior. 
There are some areas near reconstruction sites with high errors in MTCO, particularly in the northeast. 
This could reflect the fact that vegetation towards the cold and dry end of the winter temperature gradient is less sensitive to temperature change than in the Mediterranean region. 
However, the high median error for MTCO overall shows that there need to be large changes in MTCO from the prior to match the observations.

\section{Discussion}
Our final temperature reconstructions show good coherence spatially, plausible seasonal relationships, and no systematic discrepancies from pollen-based reconstructions at individual sites. 
However, the reconstructions of moisture variables, MAP and MI, are wetter than indicated by the pollen-based reconstructions. 
This was expected and is realistic. 
The atmospheric CO$_2$ concentration, [CO$_2$], was considerably lower during the LGM than it is today (180 ppm compared to 280 ppm in the PI simulations, and ca. 400 ppm today). 
Low [CO$_2$] decreases the water-use efficiency of plants and favours drought-adapted plants at the expense of trees, even without a change in climate \citep{jolly1997effect, prentice2009ecosystem}. 
Although there are methods of accounting for this direct [CO$_2$] effect \citep{prentice2017reconstructing}, statistical techniques based strictly on the application of modern analogues do not account for this impact. 
All of the pollen-based reconstructions for southern Europe from the \citet{bartlein2011pollen} data set are based on statistical reconstruction techniques. 
Application of the theoretically-based correction factor derived by \citet{prentice2017reconstructing} to the reconstructed moisture variables would be a useful next step to improve their realism.

Sites suitable for obtaining pollen records are not uniformly distributed geographically, and in any case the actual sampling of potential environments is extremely uneven in many regions of the world (Figure \ref{fig:comb_paper_ana}; \citealp{bartlein2011pollen, harrison2016what}). 
We have shown that the condition number can be used to identify an appropriate scale for interpolating the site-based data spatially, and that a scale of 400-500km appears to be appropriate for southern Europe at the LGM given the data currently available. 
This spatial scale is not uniformly appropriate, however.
The standard deviation of the reconstructions (Fig. \ref{fig:comb_paper_sd}) provides a measure of how reliable the interpolation is. 
More importantly, the standard deviation of the reconstruction could be used to determine when the interpolated values provide a realistic measure of the actual climate and when they do not. 
Establishing an acceptable threshold value for reliability would be a useful step in the creation of the kind of palaeoclimate reanalysis we are proposing here.

Whilst the values of both scales, $L_s$ and $L_t$, have been shown to be appropriate for the example shown in this paper, they are somewhat subjective.
The spatial scale, $L_s$, is chosen to give high utilisation of sparse observation data and is shown, by the condition number in Fig. \ref{fig:comb_paper_cond_num}, not to lead to a numerically inaccurate solution.
A value for $L_t$ is determined by plotting the resolution matrix for multiple $L_t$, as shown in Fig. \ref{fig:res_mat_comb_paper}; however, this only provides a range of possible values.
A more objective method for selecting $L_t$ could be developed by selecting the $L_t$ which gives the resolution matrix closest to the identity.

\section{Conclusions}
\label{sec:conclusions}
In this paper we have demonstrated a novel method for reconstructing spatially explicit palaeoclimate reconstructions from site-based data. 
The method allows the effects of each site in the dataset to be tuned by imposing a structure on the error of the prior that creates reconstructions that are spatially smooth and hence more realistic. 
By assuming that the error in the prior with respect to temperature has a given correlation month by month, it also allows the generation of a solution that is temporally smooth. 
We show that a length scale $L_t$ of 1 provides a smooth solution for the seasonal cycle, both using single sites and over multiple grid cells. 
Our analyses suggest that a spatial length scale ($L_s$) of 400km is reasonable for southern Europe at the LGM; although this is larger than the assumed source area of most of the reconstruction sites, it reflects the large-scale coherence of the regional climate change between LGM and present. 
Additional work could help to determine a more objective way to determine these length scales, but nevertheless the final climate maps appear plausible and suggest that the application of this new method should yield more robust data sets for climate-model evaluation.

\appendix
\section{Non-dimensionalisation}
\label{sec:non-dimensionalisation}
Most of the variables from the site-based reconstructions and PMIP3 have a dimension.
This can cause a problem when computing the cost function as different variables can be at different scales and it is difficult to compare different scales together computationally.
To avoid this problem we non-dimensionalise all the variables involved before computing the cost function and then re-dimensionalise the variables when the analysis has been found.

We non-dimensionalise the observation space using
\begin{linenomath*}
    \begin{equation}
        \label{eq:obs_non_dim}
        D_y(\bo{y}_i) 
        = 
        \left(
        \begin{array}{c}
            \alpha \\ D_P(P) \\ \frac{MAT}{T_s} \\ \frac{MTWA}{T_s} \\ \frac{MTCO}{T_s} \\ \frac{GDD5}{N_y T_s}
        \end{array}
        \right)
    \end{equation}
\end{linenomath*}
where $N_y$ is the number of days in a year, $T_s$ is a temperature scaling value ($5^{\circ}C$).
The function $D_P$ is defined as 
\begin{linenomath*}
    \begin{equation}
        D_P(P) = 
        \begin{cases}
            \ln \left( \frac{P \lambda}{I_{sc}} \right) + 1 & P < \frac{I_{sc}}{\lambda} \\\
            \frac{P \lambda}{I_{sc}} & \text{else}
        \end{cases}
    \end{equation}
\end{linenomath*}
where $I_{sc}$ is the solar constant ($1360.8Wm^{-2}$) and $\lambda$ is the latent heat of vaporisation of water ($2.45 MJ kg^{-1}$).
$D_P$ ensures that the method never creates a situation where $P<0$. 
Similar to the observation space, we also non-dimensionalise the state space using
\protect
\begin{linenomath*}
    \begin{equation}
        \label{eq:state_non_dim}
        D_x(\bo{x}_j) = 
        \left(
        \begin{array}{c}
            D_P(P) \\
            \frac{1}{T_s} \bo{T} \\
        \end{array}
        \right).
    \end{equation}
\end{linenomath*}

%
%
%
%
%
%
%
%

\acknowledgments
SFC was supported by a UK Natural Environment Research Programme (NERC) scholarship as part of the SCENARIO Doctoral Training Partnership at the University of Reading. 
SPH acknowledges support from the ERC-funded project GC 2.0 (Global Change 2.0: Unlocking the past for a clearer future, grant number 694481). 
ICP acknowledges support from the ERC under the European Union's Horizon 2020 research and innovation programme (grant agreement No: 787203 REALM).
This research is a contribution to the AXA Chair Programme in Biosphere and Climate Impacts and the Imperial College initiative on Grand Challenges in Ecosystems and the Environment (ICP).
NKN is supported in part by the NERC National Center for Earth Observation (NCEO).
We thank PMIP colleagues who contributed to the production of the palaeoclimate reconstructions.  
We also acknowledge the World Climate Research Programme’s Working Group on Coupled Modelling, which is responsible for CMIP, and the climate modeling groups in the Paleoclimate Modelling Intercomparison Project (PMIP) for producing and making available their model output. 
For CMIP, the U.S. Department of Energy’s Program for Climate Model Diagnosis and Intercomparison provides coordinating support and led development of software infrastructure in partnership with the Global Organization for Earth System Science Portals. The analyses and figures are based on data archived at CMIP on 12/09/18. 
We thank the Next-Generation Vegetation Modelling group for providing model code for the calculation of bioclimatic variables and for discussion of the results.

\subsection*{Code and data availability}
The reconstructions of southern Europe in this paper, as well as the code used to perform the data assimilation, are currently being archived at the University of Reading Research Data Archive (https://researchdata.reading.ac.uk/).
The CRU CL v2.0 dataset was downloaded from the University of East Anglia Climatic Research Unit and was published as \cite{new2002high}.
The PMIP3 LGM simulations \citep{braconnot2012evaluation} are available at CMIP5 archives and, for this paper, were downloaded from the Earth System Grid Federation at the Pierre Simon Laplace Institute (https://esgf-node.ipsl.upmc.fr/projects/esgf-ipsl/).
The pollen reconstructions used are available from \cite{bartlein2011pollen}.


\bibliography{library_full}

\end{document}